\documentclass[12pt]{article}
\usepackage{amsmath,amscd,amsthm,a4,amssymb}

\usepackage[russian]{babel}

\begin{document}
\begin{center}
{\Large\textbf{Нелокальная  комбинированная задача типа Бицадзе--Самарского  и  Самарского--Ионкина для  системы  псевдопараболических уравнений}}

\

И.Г. Мамедов\footnote{Институт Кибернетики им. А.И.Гусейнова НАН Азербайджана, Баку, Азербайджан.
Е-mail: ilgar-mammadov@rambler.ru
}
\end{center}

\begin{abstract}

{\it В работе рассматривается задача с условиями Бицадзе--Самарского  и  Самарского--Ионкина для  системы  псевдопараболических уравнений четвертого порядка
с разрывными коэффициентами.}

\

{\bf Ключевые слова:} псевдопараболическое уравнение, нелокальная задача.
\end{abstract}

{\bf Введение.}
В данной статье  в одном  анизотропном  пространстве С.Л. Соболева исследуется нелокальная комбинированная задача  типа Бицадзе--Самарского и  Самарского--Ионкина для одного векторного уравнения псевдопараболического типа с  доминирующей производной четвертого  порядка с негладкими коэффициентами.
Особо нужно отметить, что нелокальные задачи типа Бицадзе--Самарского  и  Самарского--Ионкина начались с работы [1] и развивались в работах [2-5] и др. Для достаточно адекватного описания большинства реальных процессов, происходщих в природе, технике и т.д. привлекаются псевдопараболические уравнения [6-8].
Многие процессы возникающие в теории фильтрации жидкости в трещиноватых средах [9], описываются псевдопараболическими уравнениями с разрывными коэффициентами. При этом математическая модель процесса дополняется нелокальными краевыми условиями типа Бицадзе--Самарского  и  Самарского--Ионкина [10]. Актуальность исследований, проводимых в этой области, объясняется появлением локальных и нелокальных граничных задач для уравнений с разрывными коэффициентами, связанных с различными прикладными задачами. Задачи такого типа возникают, в частности, при исследовании вопросов фильтрации жидкости в пористых средах, влагопереноса в грунтах, распространения импульсных лучевых волн, в различных биологических процессах и в теории обратных задач. Поэтому тематика данной работы весьма актуальна для решения многих теоретических и практических задач. Ниже для исследований таких задач предложена методика, которая исспользует современные методы теории функций и функционального анализа. Она изложена в основном применительно к псевдопараболическим уравнениям четвертого порядка с разрывными коэффициентами. При этом важным принципиальным моментам является то, что рассматриваемое  уравнение обладает, разрывными коэффициентами которые удовлетворяют только некоторым условиям типа $P$-интегрируемости и ограниченности т.е. рассмотренный псевдопараболический дифференциальный оператор не имеет традиционного сопряженного оператора. Даже в частности, например, при условиях Гурса функция Римана для такого  уравнения не может быть исследована классическим методом характеристик.

{\bf 1. Постановка задачи.}
Пусть задано уравнение
\begin{equation} \label{GrindEQ__1_}
\left(V_{1,3} u\right)\left(t,x\right)\equiv D_{t} D_{x}^{3} u \left(t,x\right)+
\underset{i+j<4}{\sum\limits _{i=0}^{1}\sum\limits _{j=0}^{3}}
\left(D_{t}^{i} D_{x}^{j} u\left(t,x\right)\right)A_{i,j} \left(t,x\right)=
$$
$$
=  Z_{1,3} \left(t,x\right),\, \left(t,x\right)\in G=
  {\left(t_{0} ,\, t_{1} \right)\times \left(x_{0} ,\, x_{1} \right)}
\end{equation}
начальное  условие
\begin{equation} \label{GrindEQ__2_}
\left(l_{0} u\right)\left(x\right)\equiv u\left(t_{0} ,x\right)=Z_{0} \left(x\right),\, \, x\in \left(x_{0} ,\, x_{1} \right),
\end{equation}
и граничные условия  типа  Бицадзе--Самарского и  Самарского--Ионкина [11]
\[
\left(l_{1}^{0} u\right)\left(t\right)\equiv u\left(t,\, x_{0} \right)\, \alpha _{1,1} +u_{x} \left(t,\, x_{0} \right)\alpha _{1,2} +u_{x\, x} \left(t,\, x_{0} \right)\alpha _{1,3} +u\left(t,\, x_{1} \right)\beta _{1,1} +
 \]
 \[
+u_{x} {\left(t,\, x_{1} \right)\beta _{1,2} +u_{xx} \left(t,\, x_{1} \right)\beta _{1,3} =\psi _{1} \left(t\right)},
\]
\begin{equation} \label{GrindEQ__3_}
{\left(l_{2}^{0} u\right)\left(t\right)\equiv u\left(t,\, x_{0} \right)\, \alpha _{2,1} +u_{x} \left(t,\, x_{0} \right)\alpha _{2,2} +u_{x\, x} \left(t,\, x_{0} \right)\alpha _{2,3} +u\left(t,\, x_{1} \right)\beta _{2,1} + }
$$
$$
+u_{x}{\left(t,\, x_{1} \right)\beta _{2,2} +u_{xx} \left(t,\, x_{1} \right)\beta _{2,3} =\psi _{2} \left(t\right)},
\end{equation}
\[
{\left(l_{3}^{0} u\right)\left(t\right)\equiv u\left(t,\, x_{0} \right)\, \alpha _{3,1} +u_{x} \left(t,\, x_{0} \right)\alpha _{3,2} +u_{x\, x} \left(t,\, x_{0} \right)\alpha _{3,3} +u\left(t,\, x_{1} \right)\beta _{3,1} + }
\]
\[
+u_{x}{\left(t,\, x_{1} \right)\beta _{3,2} +u_{xx} \left(t,\, x_{1} \right)\beta _{3,3} =\psi _{3} \left(t\right),\, \, \, t\in \left(t_{0} ,\, t_{1} \right)}.
\]

Некоторые классы граничных задач для уравнения (1) в определенном смысле ставятся аналогично известным граничным задачам для параболического уравнения $D_t u(t,x)= D^2_x u(t,x).$ Поэтому многие авторы уравнение вида (1) называют псевдопараболическим. Заметим, что рассматриваемое псевдопараболическое уравнение -- это обобщение многих модельных уравнений некоторых процессов (например, обобщенного уравнения влагопереноса, уравнения теплопроводности, уравнения колебания струны, и.т.д.).

Здесь: $u\left(t,x\right)=\left(u_{1} \left(t,x\right),...,u_{n} \left(t,x\right)\right)-n$- мерная  искомая вектор функция; $\alpha _{i,j} $ и $\beta _{i,j} $- заданные $n\times n$- мерные постоянные  матрицы, $A_{i,j} \left(t,x\right)$- измеримые на $G$ матричные функции  порядка $n\times n$,  удовлетворяющие условиям:

$A_{0,j} \left(t,x\right)\in L_{p} \left(G\right),$ $j=0,1,2$ и существуют функции $A_{1,j}^{0} \left(x\right)\in L_{p} \left(x_{0} ,\, x_{1} \right)$ и $A_{0,3}^{0} \left(t\right)\in L_{p} \left(t_{0} ,\, t_{1} \right)$ такие что,  выполнены условия
\[
\left\| A_{1,j} \left(t,\, x\right)\right\| \le A_{1,j}^{0} \left(x\right),\, \, j=0,1,2
\]
и
\[
\left\| A_{0,3} \left(t,\, x\right)\right\| \le A_{0,3}^{0} \left(t\right)
\]
почти   всюду на $G$, где $\left\| \ \ \right\| $ евклидова норма соответствующей  матрицы (или вектора); $Z_{0} \left(x\right)\in W_{p,n}^{(3)} \left(x_{0} ,\, x_{1} \right)$, а также  $\psi _{1} \left(t\right),\, \psi _{2} \left(t\right),\, \psi _{3} \left(t\right)\, \in W_{p,n}^{\left(1\right)} \left(t_{0} ,\, t_{1} \right)$- заданные $n$- мерные вектор -функции, где  $W_{p,n}^{\left(m\right)} \left(y_{0} ,\, y_{1} \right)$- пространство  $n$- мерных вектор-функций $Z\left(y\right)=\left(Z_{1} \left(y\right),...,Z_{n} \left(y\right)\right)$, имеющих в смысле С.Л. Соболева производные  $Z'\left(y\right),...,Z^{\left(m\right)} \left(y\right)\in L_{p,n} \left(y_{0,\, } y_{1} \right)$, а  $L_{p,n} \left(y_{0} ,\, y_{1} \right)$ пространство  всех строчных векторов  $Z\left(y\right)=\left(Z_{1} \left(y\right),...,Z_{n} \left(y\right)\right)$ с элементами из  $Z_{i} \left(y\right)\in L_{p} \left(y_{0} ,y_{1} \right),\, i=1,...,n$. Решение задачи (1)-(3) будем искать в  пространстве С.Л.Соболева
\[
W_{p,n}^{\left(1,3\right)} \left(G\right)=\left\{u\in L_{p,n} \left(G\right)\backslash \, \, D_{t}^{i} \, D_{x}^{j} \, \, u\in L_{p,n} \left(G\right),\, \, i=0,1,\, \, \, j=0,1,2,3\right\}
\]
с доминирующей производной  $D_{t} \, D_{x}^{3} $. Норму в нем определим равенством
\[
\left\|u\right\| _{W_{p,n}^{\left(1,3\right)} \left(G\right)} =\sum _{i=0}^{1}\sum _{j=0}^{3}\left\| D_{t}^{i}
D_{x}^{j} \, u\right\| _{L_{p,n} \left(G\right)}.
\]

Очевидно, что правые части
$Z_{0} \left(x\right)$ и $\psi _{i} \left(t\right),\, i=1,2,3$  условий (2) и (3) должны удовлетворять условиям согласования
\begin{equation} \label{GrindEQ__4_}
\left\{\begin{array}{l}
Z_{0} \left(x_{0} \right)\alpha _{1,1} +Z'_{0} \left(x_{0} \right)\alpha _{1,2} +Z''_{0} \left(x_{0} \right)\alpha _{1,3} +Z_{0} \left(x_{1} \right)\beta _{1,1} +\\ +Z'_{0} \left(x_{1} \right)\beta _{1,2} +Z''_{0} \left(x_{1} \right)\beta _{1,3} = \psi _{1} \left(t_{0} \right) \\
Z_{0} \left(x_{0} \right)\alpha _{2,1} +Z'_{0} \left(x_{0} \right)\alpha _{2,2} +Z''_{0} \left(x_{0} \right)\alpha _{2,3} +Z_{0} \left(x_{1} \right)\beta _{2,1} +\\
+Z'_{0} \left(x_{1} \right)\beta _{2,2} +Z''_{0} \left(x_{1} \right)\beta _{2,3} =\psi _{2} \left(t_{0} \right)
\\
Z_{0} \left(x_{0} \right)\alpha _{3,1} +Z'_{0} \left(x_{0} \right)\alpha _{3,2} +Z''_{0} \left(x_{0} \right)\alpha _{3,3} +Z_{0} \left(x_{1} \right)\beta _{3,1} +\\
+Z'_{0} \left(x_{1} \right)\beta _{3,2} +Z''_{0} \left(x_{1} \right)\beta _{3,3}  =\psi _{3} \left(t_{0} \right) \end{array}\right.  \end{equation}

Наличие  условий согласования  означает, что условиями (2) и (3)  заданы также некоторые  излишние информации о решении. Поэтому желательно, чтобы в  постановке этой задачи не  оказалась необходимости к  условиям  типа согласования.

Для получения таких условий мы будем дифференцировать условия (3) по $t$. Тогда  получим
\[
{\left(l_{1} \, u\right)\left(t\right)\equiv u_{t} \left(t,\, x_{0} \right)\alpha _{1,1} +u_{tx} \left(t,\, x_{0} \right)\alpha _{1,2} +u_{txx} \left(t,\, x_{0} \right)\alpha _{1,3} +u_{t} \left(t,\, x_{1} \right)\beta _{1,1} +}
\]
\[
{+u_{tx} \left(t,\, x_{1} \right)\beta _{1,2} +u_{txx} \left(t,x_{1} \right)\beta _{1,3} =Z_{1} \left(t\right)=\psi '_{1} \left(t\right)}
\]
\begin{equation} \label{GrindEQ__5_}
 {\left(l_{2} \, u\right)\left(t\right)\equiv u_{t} \left(t,\, x_{0} \right)\alpha _{2,1} +u_{tx} \left(t,\, x_{0} \right)\alpha _{2,2} +u_{txx} \left(t,\, x_{0} \right)\alpha _{2,3} +u_{t} \left(t,\, x_{1} \right)\beta _{2,1} +}
 $$
 $$
 {+u_{tx} \left(t,\, x_{1} \right)\beta _{2,2} +u_{txx} \left(t,x_{1} \right)\beta _{2,3} =Z_{2} \left(t\right)=\psi '_{2} \left(t\right)}
\end{equation}
\[
{\left(l_{3} \, u\right)\left(t\right)\equiv u_{t} \left(t,\, x_{0} \right)\alpha _{3,1} +u_{tx} \left(t,\, x_{0} \right)\alpha _{3,2} +u_{txx} \left(t,\, x_{0} \right)\alpha _{3,3} +u_{t} \left(t,\, x_{1} \right)\beta _{3,1} +}
\]
\[
{+u_{tx} \left(t,\, x_{1} \right)\beta _{3,2} +u_{txx} \left(t,x_{1} \right)\beta _{3,3} =Z_{3} \left(t\right)=\psi '_{3} \left(t\right)}
\]

Причем здесь будем требовать, чтобы  выполнялись  лишь условия
\[
Z_{i} \in L_{p,n} \left(t_{0} ,t_{1} \right),\, \, i=1,2,3.
\]

Очевидно, что если  $u\in W_{p,n}^{\left(1,3\right)} \left(G\right)$ есть решение задачи (1), (2), (5), то она является также  решением задачи (1)-(3)  при
\begin{equation} \label{GrindEQ__6_}
 \psi _{1} \left(t\right)=\int\limits _{t_{0} }^{t}Z_{1} \left(\tau \right)d\tau  +Z_{0} \left(x_{0} \right)\alpha _{1,1} +Z'_{0} \left(x_{0} \right)\alpha _{1,2} +
 $$$$
  +Z''_{0} \left(x_{0} \right)\alpha _{1,3} +Z_{0} \left(x_{1} \right)\beta _{1,1} +Z'_{0} \left(x_{1} \right)\beta _{1,2}
  {+Z''_{0} \left(x_{1} \right)\beta _{1,3} ;}
 $$
 $$
   {\psi _{2} \left(t\right)=\int\limits _{t_{0} }^{t}Z_{2} \left(\tau \right)d\tau  +Z_{0} \left(x_{0} \right)\alpha _{2,1} +Z'_{0} \left(x_{0} \right)\alpha _{2,2} +Z''_{0} \left(x_{0} \right)\alpha _{2,3} +}
 $$
 $$
    {+Z_{0} \left(x_{1} \right)\beta _{2,1} +Z'_{0} \left(x_{1} \right)\beta _{2,2} +Z''_{0} \left(x_{1} \right)\beta _{2,3} ;}
 $$
 $$
     {\psi _{3} \left(t\right)=\int\limits _{t_{0} }^{t}Z_{3} \left(\tau \right)d\tau  +Z_{0} \left(x_{0} \right)\alpha _{3,1} +Z'_{0} \left(x_{0} \right)\alpha _{3,2} +Z''_{0} \left(x_{0} \right)\alpha _{3,3} +}
$$
$$
+Z_{0} \left(x_{1} \right)\beta _{3,1} +Z'_{0} \left(x_{1} \right)\beta _{3,2} {+Z''_{0} \left(x_{1} \right)\beta _{3,3}}.
 \end{equation}

Отметим, что для  функций (6) условия согласования (4) выполняются тривиальным  образом. Верно и обратное, т.е. если  $u\in W_{p,n}^{\left(1,3\right)} \left(G\right)$ является  решением задачи  (1)-(3), то она  является также решением задачи  (1), (2), (5) при  $Z_{1} \left(t\right)=\psi '_{1} \left(t\right),\, \, Z_{2} \left(t\right)=\psi '_{2} \left(t\right)$ и $Z_{3} \left(t\right)=\psi '_{3} \left(t\right)$. Иначе  говоря,  задачи (1)-(3) и (1), (2), (5) эквивалентны в $W_{p,n}^{\left(1,3\right)} $. При этом  решение задачи (1), (2),(5) есть  решение  задачи (1)-(3) при некоторых  $\psi _{1} \left(t\right),\, \psi _{2} \left(t\right),\, \psi _{3} \left(t\right)$ удовлетворяющих условиям согласования  автоматическим образом.

Однако, задача (1), (2), (5) по постановке более естественна,  чем (1)-(3). Поэтому в дальнейшем  будем исследовать только  задачу  (1), (2),(5).\\

\begin{center}
{\bf 2. Интегральное  представление функции в пространстве\\ С.Л.Соболева  посредством определяющих операторов  при\\ исследовании нелокальной задачи}
\end{center}

Различные классы нелокальных задач типа Бицадзе--Самарского для дифференциальных уравнений являлись результатом исследований многих математиков [12-15] и др. Кроме того, нужно отметить, что в литературе исследованы разнообразные нелокальные задачи типа Самарского--Ионкина для различных классов дифференциальных уравнений [16-17]. Такие нелокальные краевые задачи типа Бицадзе--Самарского  и  Самарского--Ионкина в модифицированных трактовках рассматривались в работах автора [18-19].

Рассмотрим задачу (1), (2),(5). Эту задачу напишем в операторном виде
\begin{equation} \label{GrindEQ__7_}
Vu=Z,
\end{equation}
Здесь:
\[
{ V}=\left({\rm V}_{{\rm 1,3}} ,\, l_{0} ,\, l_{1} ,\, l_{2} ,\, l_{3} \right),
\]
\[
Z=\left(Z_{1,3} \left(t,x\right),\, Z_{0} \left(x\right),\, Z_{1} \left(t\right),\, Z_{2} \left(t\right),\, Z_{3} \left(t\right)\right)\in E_{p,n}^{\left(1,3\right)} ,
\]
\[
E_{p,n}^{\left(1,3\right)} \equiv L_{p,n} \left(G\right)\times W_{p,n}^{\left(3\right)} \left(x_{0} ,\, x_{1} \right)\times L_{p,n} \left(t_{0} ,t_{1} \right)\times L_{p,n} \left(t_{0} ,t_{1} \right)\times L_{p,n} \left(t_{0} ,t_{1} \right).
\]

Норму в пространстве $E_{p,n}^{\left(1,3\right)} $ будем считать определенной  естественным образом при  помощи  равенства
\[
\left\| Z\right\| _{E_{p,n}^{\left(1,3\right)} } =\left\| Z_{1,3} \right\| _{L_{p,n} \left(G\right)} +\left\| Z_{0} \right\| _{W_{p,n}^{\left(3\right)} \left(x_{0} ,\, x_{1} \right)} +\left\| Z_{1} \right\| _{L_{p,n} \left(t_{0} ,t_{1} \right)}+
\]
\[
+ \left\| Z_{2} \right\| _{L_{p,n} \left(t_{0} ,\, t_{1}\right)}+\left\| Z_{3} \right\| _{L_{p,n} \left(t_{0} ,\, t_{1} \right)} .
\]

Можно доказать, что оператор  $V:W_{p,n}^{\left(1,3\right)} \left(G\right)\to E_{p,n}^{\left(1,3\right)} $ линеен и  ограничен.

\textbf{Определение.} {\it Если для любого  $Z\in E_{p,n}^{\left(1,3\right)} $ уравнение (7)  (задача (1), (2), (5)) имеет единственное  решение $u\in W_{p,n}^{\left(1,3\right)} \left(G\right)$ такое, что
\begin{equation}\label{GrindEQ__8_}
\left\| u\right\| _{W_{p,n}^{\left(1,3\right)} \left(G\right)} \le M\left\| Z\right\| _{E_{p,n}^{\left(1,3\right)} },
\end{equation}
то будем говорить, что  оператор $V$ уравнения (7) (задачи (1), (2), (5)) осуществляет гомеоморфизм  между $W_{p,n}^{\left(1,3\right)} \left(G\right)$ и $E_{p,n}^{\left(1,3\right)} $, или  будем говорить, что задача (1), (2), (5) везде корректно  разрешима, где $M>0$ постоянное (независящее от коэффициентов  краевой задачи (1), (2), (5), а также от размерностей области $G$), независящее от $Z\in E_{p,n}^{\left(1,3\right)} $.}

В случае, когда $V$ есть гомеоморфизм между $W_{p,n}^{\left(1,3\right)} \left(G\right)$ и $E_{p,n}^{\left(1,3\right)} $, оператор $V$ обладает ограниченной обратной  определенной  на $E_{p,n}^{\left(1,3\right)} $.

В современной теории дифференциальных уравнений особое значение имеет вопрос о выявлении классов задач, операторы которых осуществляют гомеоморфизм между определенными парами банаховых пространств. Такие гомеоморфизмы выявлены в работах Ю.М.Березанского и Я.А.Ройтберга [20], Н.В.Житарашу [21], С.С.Ахиева [22], автора [23-24] и др. для некоторых классов дифференциальных уравнений с частными производными.

Задачу (1), (2), (5) будем исследовать  при помощи  интегральных  представлений специального вида для функций $u\in W_{p,n}^{\left(1,3\right)} \left(G\right)$. Для  функций $u\in W_{p,n}^{\left(1,3\right)} \left(G\right)$ можно найти различные интегральные  представления. Функцию $u\in W_{p,n}^{\left(1,3\right)} \left(G\right)$ можно представить, например,  в виде
\begin{equation} \label{GrindEQ__9_}
u(t,x)=u(t_0,x)+\int\limits\limits_{t_0}^t\left[u_t(\tau,x_0)+(x-x_0)u_{tx}(\tau,x_0)+\dfrac{(x-x_0)^2}{2!}u_{txx}(\tau,x_0)\right]d\tau+
$$
$$
+\int\limits\limits_{t_0}^t\int\limits\limits_{x_0}^xu_{txxx}(\tau,\xi)\dfrac{(x-\xi)^2}{2!}d\tau d\xi,
\end{equation}
посредством следов  $u\left(t_{0} ,x\right),\, \, u_{t} \left(t,x_{0} \right),\, \, u_{tx} \left(t,\, x_{0} \right),\, \, u_{txx} \left(t,x_{0}\right)$ и доминирующей  производной  $u_{txxx} \left(t,x\right)$. Существуют также  представления других видов, позволяющих  выразить функцию $u\left(t,x\right)$ посредством некоторых следов и  частных производных. Однако,  мы будем ставить такой  вопрос. Нельзя ли для функций $u\in W_{p,n}^{\left(1,3\right)} \left(G\right)$ найти представление, которое позволило бы определить  функцию $u\left(t,x\right)$ однозначно по  значениям $\left(B_{k} \, u\right)$, $k=1,...,m,$ некоторых заданных операторов  $B_{k} ,\, k=1,..,m,$ на этой функции  $u\left(t,x\right)$. Если, например, функцию $u\in W_{p,n}^{\left(1,3\right)} \left(G\right)$ можно  было  восстановить однозначно по  значениям
$\left(V_{1,3} u\right)\left(t,x\right),\, \, \left(l_{0} \, u\right)\left(x\right)$, $\left(l_{1} u\right)\left(t\right),\, \, \left(l_{2} u\right)\left(t\right)$ и $\left(l_{3} u\right)\left(t\right)$, операторов $V_{1,3} ,\, l_{0} ,\, l_{1} ,\, l_{2} $ и $l_{3} $ на  этой функции, то нахождение решения задачи (1), (2), (5), не  представляло бы собою особую  трудность. Однако, нахождение  подобного представления равносильно построению обратного  оператора для оператора $V$ этой  задачи. Поэтому этот вопрос  мы  будем ставить в следующем  ослабленном виде. Найти для  функции $u\in W_{p,n}^{\left(1,3\right)} \left(G\right)$  представление, по которому можно  было однозначно восстановить  функцию  $u\left(t,x\right)$ по значениям  $\left(l_{0} u\right)\left(x\right),$ $\left(l_{1} u\right)\left(t\right),$ $\left(l_{2} u\right)\left(t\right),$ $\left(l_{3} u\right)\left(t\right)$ и $D_{t} D_{x}^{3}$ $ u\left(t,x\right)$.

 Для изучения этого вопроса мы будем использовать формулу (9). Сначала  формулу (9) запишем в  виде
\[
 u\left(t,x\right)=\left(Qb\right)\left(t,x\right)\equiv \varphi \left(x\right)+
 \]
 \begin{equation} \label{GrindEQ__10_}
 +\int\limits _{t_{0} }^{t}\left[b_{1,0} \left(\tau \right)
 +\left(x-x_{0} \right)b_{1,1} \left(\tau \right)+\frac{\left(x-x_{0} \right)^{2} }{2!} b_{1,2} \left(\tau \right)\right] d\tau +
\end{equation}
 \[
+ \int\limits _{t_{0} }^{t}\int\limits _{x_{0} }^{x}\frac{\left(x-\xi \right)^{2} }{2!} b_{1,3} \left(\tau ,\xi \right)d\tau \, d\xi ,\, \left(t,x\right)\in G  ,
\]
где
\[
b=\left(b_{1,3} \left(t,x\right),\, \varphi \left(x\right),\, b_{1,0} \left(t\right),\, b_{1,1} \left(t\right),\, b_{1,2} \left(t\right)\right)
\]
и
\begin{equation} \label{GrindEQ__11_}
\begin{array}{c} {\varphi \left(x\right)=u\left(t_{0} ,x\right)}, \\ {b_{1,0} \left(t\right)=u_t\left(t,x_{0} \right)}, \\ {b_{1,1} \left(t\right)=u_{tx}\left(t,x_{0} \right)}, \\ {b_{1,2} \left(t\right)=u_{txx}\left(t,x_{0} \right)}, \\ {b_{1,3} \left(t\right)=u_{txxx}\left(t,x\right)} . \end{array}
\end{equation}

Известно,  что для  любой функции $u\in W_{p,n}^{\left(1,3\right)} \left(G\right)$ существует  единственная пятерка
\[
b=\left(b_{1,3} ,\, \varphi ,\, b_{1,0} ,\, b_{1,1} ,\, b_{1,2} \right)\in E_{p,n}^{\left(1,3\right)} ,
\]
посредством которой эту  функцию можно представить в  виде (10).

Можно доказать, что верно и обратное. Иначе  говоря, для  любой заданной пятерки
$b=\left(b_{1,3} ,\, \varphi ,\, b_{1,0} ,\, b_{1,1} ,\, b_{1,2} \right)$ из $E_{p,n}^{\left(1,3\right)} $  функция $u\left(t,x\right)$,  определяемая  равенством (10)  принадлежит  пространству $W_{p,n}^{\left(1,3\right)} \left(G\right)$ и обладает следами $u\left(t_{0} ,x\right)$, $u_{t} \left(t,x_{0} \right),\, u_{tx}\left(t,\, x_{0} \right)$,
$u_{txx} \left(t,\, x_{0} \right)$, а также  доминирующей  производной $u_{txxx} \left(t,\, x\right)$ определяемой  равенствами (11).

Боле того, при этом можно найти положительные  постоянные $C_{1} $ и  $C_{2} $ такие, что
\begin{equation} \label{GrindEQ__12_}
C_{1} \left\| b\right\| _{E_{p,n}^{\left(1,3\right)} } \le \left\| Qb\right\| _{W_{p,n}^{\left(1,3\right)} \left(G\right)} \le C_{2} \left\| b\right\| _{E_{p,n}^{\left(1,3\right)} } ,\, \, \, \, \forall b\in E_{p,n}^{\left(1,3\right)}.
\end{equation}

Иначе говоря, оператор $Q$ осуществляет гомеоморфизм  между пространствами
$E_{p,n}^{\left(1,3\right)} $ и $W_{p,n}^{\left(1,3\right)} \left(G\right)$. Поэтому $W_{p,n}^{\left(1,3\right)} \left(G\right)=E_{p,n}^{\left(1,3\right)} $ в смысле гомеоморфизма.

Теперь элемент $b=(b_{1,3},\varphi,b_{1,0},b_{1,1},b_{1,2})\in E_{p,n}^{(1,3)}$ постараемся  выбирать  таким образом, чтобы   соответствующая функция (10)  удовлетворяла условиям (2) и (5). Для этого подставим  функцию (10) в условиях (2) и (5).

Тогда из (2) получим, что
\begin{equation}\label{GrindEQ__13_}
\left(l_{0} \, u\right)\left(x\right)=\varphi \left(x\right),
\end{equation}

Из условий (5) получим  следующие условия:
\[
 \left(l_{1}  u\right)\left(t\right)=b_{1,0} \left(t\right)\alpha _{1,1} +b_{1,1} \left(t\right)\alpha _{1,2} +b_{1,2} \left(t\right)\alpha _{1,3} +
 \]
 \[
 +\left[b_{1,0} \left(t\right)+\Delta b_{1,1} \left(t\right)+\frac{\Delta ^{2} }{2!}
 b_{1,2} \left(t\right)+ \int\limits _{x_{0} }^{x_{1} }\frac{\left(x_{1} -\xi \right)^{2} }{2!} b_{1,3} \left(t, \xi \right)d\xi  \right]\beta _{1,1} +
 \]
 \[
 +\left[b_{1,1} \left(t\right)+\Delta b_{1,2} \left(t\right)+\int\limits _{x_{0} }^{x_{1} }\left(x_{1} -\xi \right)b_{1,3} \left(t, \xi \right)d\xi  \right]
  \beta _{1,2} +
  \]
  \[
  +\left[b_{1,2} \left(t\right)+\int\limits _{x_{0} }^{x_{1} }b_{1,3} \left(t,\, \xi \right)d\xi  \right]\beta _{1,3} ;
\]
\[
\left(l_{2} \, u\right)\left(t\right)=b_{1,0} \left(t\right)\alpha _{2,1} +b_{1,1} \left(t\right)\alpha _{2,2} +b_{1,2} \left(t\right)\alpha _{2,3} +
\]
\[
+\left[b_{1,0} \left(t\right)+\Delta b_{1,1} \left(t\right)+\frac{\Delta ^{2} }{2!} b_{1,2} \left(t\right)+
\int\limits _{x_{0} }^{x_{1} }\frac{\left(x_{1} -\xi \right)^{2} }{2!} b_{1,3} \left(t, \xi \right)d\xi  \right]\beta _{2,1} +
\]
\[
+\left[b_{1,1} \left(t\right)+\Delta b_{1,2} \left(t\right)+\int\limits _{x_{0} }^{x_{1} }\left(x_{1} -\xi \right)b_{1,3} \left(t,\, \xi \right)d\xi  \right]
\beta _{2,2} +
\]
\[
+\left[b_{1,2} \left(t\right)+\int\limits _{x_{0} }^{x_{1} }b_{1,3} \left(t,\, \xi \right)d\xi  \right]\beta _{2,3} ;
\]
\begin{equation} \label{GrindEQ__14_}
\left(l_{3}  u\right)\left(t\right)=b_{1,0} \left(t\right)\alpha _{3,1} +b_{1,1} \left(t\right)\alpha _{3,2} +b_{1,2} \left(t\right)\alpha _{3,3}+
 $$
 $$
 +\left[b_{1,0} \left(t\right)+\Delta b_{1,1} \left(t\right)+\frac{\Delta ^{2} }{2!} b_{1,2} \left(t\right)+ \int\limits _{x_{0} }^{x_{1} }\frac{\left(x_{1} -\xi \right)^{2} }{2!} b_{1,3} \left(t,,\, \xi \right)d\xi  \right]\beta _{3,1} +
 $$
 $$
 +\left[b_{1,1} \left(t\right)+\Delta b_{1,2} \left(t\right)+\int\limits _{x_{0} }^{x_{1} }\left(x_{1} -\xi \right)b_{1,3} \left(t,\, \xi \right)d\xi  \right]\beta _{3,2}+
 $$$$
  +\left[b_{1,2} \left(t\right)+\int\limits _{x_{0} }^{x_{1} }b_{1,3} \left(t,\, \xi \right)d\xi  \right]\beta _{3,3} ,
\end{equation}
где $\Delta =x_{1} -x_{0} $. Равенства (14)   запишем в виде
\[
\left(l_{1} \, u\right)\left(t\right)=b_{1,0} \left(t\right)\left(\alpha _{1,1} +\beta _{1,1} \right)+b_{1,1} \left(t\right)\left(\alpha _{1,2} +\Delta \beta _{1,1} +\beta _{1,2} \right)+
\]
\[
+b_{1,2} \left(t\right)\left(\alpha _{1,3} +\frac{\Delta ^{2} }{2!} \beta _{1,1} +\Delta \beta _{1,2} +\beta _{1,3} \right)+
\]\[
+ \int\limits _{x_{0} }^{x_{1} }b_{1,3} \left(t,\, \xi \right) \left[\beta _{1,1} \frac{\left(x_{1} -\xi \right)^{2} }{2!} +\beta _{1,2} \left(x_{1} -\xi \right)+\beta _{1,3} \right]d\xi;
\]
\[
\left(l_{2} \, u\right)\left(t\right)=b_{1,0} \left(t\right)\left(\alpha _{2,1} +\beta _{2,1} \right)+b_{1,1} \left(t\right)\left(\alpha _{2,2} +\Delta \beta _{2,1} +\beta _{2,2} \right)+
\]
\[
+b_{1,2} \left(t\right)\left(\alpha _{2,3} +\frac{\Delta ^{2} }{2!} \beta _{2,1} +\Delta \beta _{2,2} +\beta _{2,3} \right)+
 $$
 $$
 +\int\limits _{x_{0} }^{x_{1} }b_{1,3} \left(t,\, \xi \right) \left[\beta _{2,1} \frac{\left(x_{1} -\xi \right)^{2} }{2!}  +\beta _{2,2} \left(x_{1} -\xi \right)+\beta _{2,3} \right]d\xi ;
\]
\begin{equation} \label{GrindEQ__15_}
\left(l_{3} \, u\right)\left(t\right)=b_{1,0} \left(t\right)\left(\alpha _{3,1} +\beta _{3,1} \right)+b_{1,1} \left(t\right)\left(\alpha _{3,2} +\Delta \beta _{3,1} +\beta _{3,2} \right)+
$$
$$
+b_{1,2} \left(t\right)\left(\alpha _{3,3} +\frac{\Delta ^{2} }{2!} \beta _{3,1} +\Delta \beta _{3,2} +\beta _{3,3} \right)+
$$
$$
+\int\limits _{x_{0} }^{x_{1} }b_{1,3} \left(t,\, \xi \right) \left[\beta _{3,1} \frac{\left(x_{1} -\xi \right)^{2} }{2!}  +\beta _{3,2} \left(x_{1} -\xi \right)
+\beta _{3,3}  \right] d\xi.
\end{equation}

 Введем следующие обозначения
\[
\gamma _{1,1} =\alpha _{1,1} +\beta _{1,1} ;\, \, \gamma _{1,2} =\alpha _{1,2} +\Delta\beta _{1,1} +\beta_{1,2} ;
\]
\[
\gamma _{1,3} =\alpha _{1,3} +\frac{\Delta ^{2} }{2!} \beta _{1,1} +\Delta \beta _{1,2} +\beta _{1,3} ;\, \, \, \, \, \gamma _{2,1} =\alpha _{2,1} +\beta _{2,1} ;
\]
\[
\gamma _{2,2} =\alpha _{2,2} +\Delta \beta _{2,1} +\beta _{2,2} ;
\
\gamma _{2,3} =\alpha _{2,3} +\frac{\Delta ^{2} }{2!} \beta _{2,1} +\Delta \beta _{2,2} +\beta _{2,3} ;
\]
\[
\gamma _{3,1} =\alpha _{3,1} +\beta _{3,1} ;
\
\gamma _{3,2} =\alpha _{3,2} +\Delta \beta _{3,1} +\beta _{3,2} ;
\]
\[
\gamma _{3,3} =\alpha _{3,3} +\frac{\Delta ^{2} }{2!} \beta _{3,1} +\Delta\beta _{3,2} +\beta _{3,3} ;
\]
и
\[
a_{1} \left(x\right)=\beta _{1,1} \frac{\left(x_{1} -x\right)^{2} }{2!} +\beta _{1,2} \left(x_{1} -x\right)+\beta _{1,3} ;
\]
\[
a_{2} \left(x\right)=\beta _{2,1} \frac{\left(x_{1} -x\right)^{2} }{2!} +\beta _{2,2} \left(x_{1} -x\right)+\beta _{2,3} ;
\]
\[
a_{3} \left(x\right)=\beta _{3,1} \frac{\left(x_{1} -x\right)^{2} }{2!} +\beta _{3,2} \left(x_{1} -x\right)+\beta _{3,3} ;
\]

Тогда равенства (15)  примут  вид
\[
b_{1,0} \left(t\right)\gamma _{1,1} +b_{1,1} \left(t\right)\gamma _{1,2} +b_{1,2} \left(t\right)\gamma _{1,3} =\left(l_{1} u\right)\left(t\right)-\int\limits _{x_{0} }^{x_{1} }b_{1,3} \left(t,\, \xi \right)a_{1} \left(\xi \right)d\xi  ;
\]
\[
b_{1,0} \left(t\right)\gamma _{2,1} +b_{1,1} \left(t\right)\gamma _{2,2} +b_{1,2} \left(t\right)\gamma _{2,3} =\left(l_{2} u\right)\left(t\right)-\int\limits _{x_{0} }^{x_{1} }b_{1,3} \left(t,\, \xi \right)a_{2} \left(\xi \right)d\xi  ;
\]
\begin{equation} \label{GrindEQ__16_}
b_{1,0} \left(t\right)\gamma _{3,1} +b_{1,1} \left(t\right)\gamma _{3,2} +b_{1,2} \left(t\right)\gamma _{3,3} =\left(l_{3} u\right)\left(t\right)-\int\limits _{x_{0} }^{x_{1} }b_{1,3} \left(t,\, \xi \right)a_{3} \left(\xi \right)d\xi
\end{equation}

Равенства (16) можно  рассматривать как систему линейных  алгебраических уравнений  относительно неизвестных $b_{1,0} \left(t\right)$, $b_{1,1} \left(t\right)$ и {$b_{1,2} \left(t\right)$}.

При помощи $3n$-мерной матрицы
\[
\gamma =\left(\begin{array}{l} {\gamma _{1,1} \, \, \, \, \gamma _{2,1} \, \, \, \gamma _{3,1} } \\ {\gamma _{1,2} \, \, \, \gamma _{2,2} \, \, \, \gamma _{3,2} } \\ {\gamma _{1,3} \, \, \, \, \gamma _{2,3\, \, \, } \gamma _{3,3} } \end{array}\right)
\]
систему (16) запишем в виде векторного уравнения
$$
 \left(b_{1,0} \left(t\right),\, \, b_{1,1} \left(t\right),\, b_{1,2} \left(t\right)\right)\gamma =
 \left(\left(l_1 u\right)\left(t\right)-\int\limits _{x_{0} }^{x_{1}}
 b_{1,3} \left(t,\xi \right)a_{1} \left(\xi \right)d\xi  ,\right.
$$
$$
\left. \left(l_{2} \, u\right)\left(t\right)-
 \int\limits _{x_{0} }^{x_{1} } b_{1,3} \left(t, \xi \right)a_{2} \left(\xi \right)d\xi  ,
 \left(l_{3} u\right)\left(t\right)- \int\limits _{x_{0} }^{x_{1} }b_{1,3} \left(t, \xi \right)a_{3} \left(\xi \right)d\xi  \right).
\eqno(16^*)
$$

Пусть матрица $\gamma $  обратима. Обратную матрицу $\gamma ^{-1} $  обозначим через
\[
\gamma ^{-1} =\left(\begin{array}{l} {K_{1,1} \, \, \, \, K_{1,2} \, \, \, K_{1,3} } \\ {K_{2,1} \, \, \, K_{2,2} \, \, \, K_{2,3} } \\ {K_{3,1} \, \, \, \, K_{3,2\, \, } K_{3,3} } \end{array}\right)
\]
причем  каждая из  $K_{i,j} $ является  некоторой  квадратной постоянной  матрицей порядка $n$. Тогда  из (16*) получим
\[
\left(b_{1,0} \left(t\right),b_{1,1} \left(t\right), b_{1,2} \left(t\right)\right) =
\left(\left(l_1 u\right)\left(t\right)-\int\limits _{x_{0} }^{x_{1} } b_{1,3} \left(t,\xi \right)a_{1} \left(\xi \right)d\xi, \right.
\]
\[
\left(l_{2}  u\right)\left(t\right)-\int\limits _{x_{0} }^{x_{1} }b_{1,3} \left(t,\, \xi \right)a_{2} \left(\xi \right)d\xi  , \left(l_{3} u\right)\left(t\right)-\left. \int\limits _{x_{0} }^{x_{1} }b_{1,3} \left(t,\, \xi \right)a_{3} \left(\xi \right)d\xi  \right)\times
\]
\[
\times\left(\begin{array}{l} {K_{1,1} \, \, \, \, K_{1,2} \, \, \, K_{1,3} } \\ {K_{2,1} \, \, \, K_{2,2} \, \, \, K_{2,3} } \\ {K_{3,1} \, \, \, \, K_{3,2\, \, } K_{3,3} } \end{array}\right)
\]

Поэтому
\[
b_{1,0} \left(t\right)=\left[\left(l_{1} u\right)\left(t\right)-\int\limits _{x_{0} }^{x_{1} }b_{1,3} \left(t,\, \xi \right)a_{1} \left(\xi \right)d\xi  \right]K_{1,1}+
\]
\[
 +\left[\left(l_{2} u\right)\left(t\right)-\int\limits _{x_{0} }^{x_{1} }b_{1,3} \left(t,\, \xi \right)a_{2} \left(\xi \right)d\xi  \right]K_{2,1} +
\]
\[
+\left[\left(l_{3} u\right)\left(t\right)-\int\limits _{x_{0} }^{x_{1} }b_{1,3} \left(t,\, \xi \right)a_{3} \left(\xi \right)d\xi  \right]K_{3,1} ;
\]
\[
b_{1,1} \left(t\right)=\left[\left(l_{1} u\right)\left(t\right)-\int\limits _{x_{0} }^{x_{1} }b_{1,3} \left(t,\, \xi \right)a_{1} \left(\xi \right)d\xi  \right]K_{1,2} +
\]
\[
+\left[\left(l_{2} u\right)\left(t\right)-\int\limits _{x_{0} }^{x_{1} }b_{1,3} \left(t,\, \xi \right)a_{2} \left(\xi \right)d\xi  \right]K_{2,2}+
\]
\[
+\left[\left(l_{3} u\right)\left(t\right)-\int\limits _{x_{0} }^{x_{1} }b_{1,3} \left(t,\, \xi \right)a_{3} \left(\xi \right)d\xi  \right]K_{3,2} ;
\]
\begin{equation} \label{GrindEQ__17_}
b_{1,2} \left(t\right)=\left[\left(l_{1} u\right)\left(t\right)-\int\limits _{x_{0} }^{x_{1} }b_{1,3} \left(t,\, \xi \right)a_{1} \left(\xi \right)d\xi  \right]K_{1,3} +
$$
$$
+\left[\left(l_{2} u\right)\left(t\right)-\int\limits _{x_{0} }^{x_{1} }b_{1,3} \left(t,\, \xi \right)a_{2} \left(\xi \right)d\xi  \right]K_{2,3} +
$$
$$
+\left[\left(l_{3} u\right)\left(t\right)-\int\limits _{x_{0} }^{x_{1} }b_{1,3} \left(t,\, \xi \right)a_{3} \left(\xi \right)d\xi  \right]K_{3,3} ;
\end{equation}

Теперь выражения (17) учтем в равенстве (10). Тогда получим
\begin{equation} \label{GrindEQ__18_}
u\left(t,x\right)=\left(l_{0} u\right)\left(x\right)+\int\limits _{t_{0} }^{t}\left(l_{1} \, u\right)\left(\tau \right)\left[K_{1,1} +K_{1,2} \left(x-x_{0} \right)+K_{1,3} \frac{\left(x-x_{0} \right)^{2} }{2!} \right] d\tau +
$$
$$
+\int\limits _{t_{0} }^{t}\left(l_{2} u\right)\left(\tau \right) \left[K_{2,1} +K_{2,2} \left(x-x_{0} \right)+K_{2,3} \frac{\left(x-x_{0} \right)^{2} }{2!} \right]d\tau +
$$
$$
 +\int\limits _{t_{0} }^{t}\left(l_{3} u\right)\left(\tau \right)
\left[K_{3,1} +K_{3,2} \left(x-x_{0} \right)+K_{3,3} \frac{\left(x-x_{0} \right)^{2} }{2!} \right]d\tau+
$$
$$
 +\iint \nolimits _{G}b_{1,3} \left(\tau ,\, \xi \right) \left\{E\frac{\left(x-\xi \right)^{2} }{2!} \theta \left(t-\tau \right)\theta \left(x-\xi \right)\right. -
$$
$$
-\theta \left(t-\tau \right)\Big[a_{1} \left(\xi \right)K_{1,1} +a_{2} \left(\xi \right) K_{2,1} +a_{3} \left(\xi \right)K_{3,1}+
$$
$$
 +\Big(a_{1} \left(\xi \right)K_{1,2} +a_{2} \left(\xi \right)K_{2,2} +a_{3} \left(\xi \right)K_{3,2} \Big)\left(x-x_{0} \right)+
$$
$$
+\left. \left. \Big(a_{1} \left(\xi \right)K_{1,3} +a_{2} \left(\xi \right)K_{2,3} +a_{3} \left(\xi \right)K_{3,3}
   \Big)\frac{\left(x-x_{0} \right)^{2} }{2!} \right]\right\}d\tau \, d\xi,
\end{equation}
где $\theta \left(z\right)$-функция Хэвисайда на пространстве $ \mathbb{R}$  действительных чисел, а $E$- $n\times n$- мерная единичная матрица.

Пусть
\[
\beta _{1} \left(x\right)=K_{1,1} +K_{1,2} \left(x-x_{0} \right)+K_{1,3} \frac{\left(x-x_{0} \right)^{2} }{2!} =-q_{1} \left(x\right),
\]
\[
\beta _{2} \left(x\right)=K_{2,1} +K_{2,2} \left(x-x_{0} \right)+K_{2,3} \frac{\left(x-x_{0} \right)^{2} }{2!} =-q_{2} \left(x\right),
\]
\[
\beta _{3} \left(x\right)=K_{3,1} +K_{3,2} \left(x-x_{0} \right)+K_{3,3} \frac{\left(x-x_{0} \right)^{2} }{2!} =-q_{3} \left(x\right),
\]
\begin{equation} \label{GrindEQ__19_}
R_{0} \left(\tau ,\, \xi ;\, t,x\right)=\theta \left(t-\tau \right)\times
\end{equation}
\[
\times \left[\frac{\left(x-\xi \right)^{2} }{2!} \theta \left(x-\xi \right)E+a_{1} \left(\xi \right)q_{1} \left(x\right)+a_{2} \left(\xi \right)q_{2} \left(x\right)+a_{3} \left(\xi \right)q_{3} \left(x\right) \right].
\]

Теперь формулу (18) запишем  в виде
\begin{equation} \label{GrindEQ__20_}
u\left(t,x\right)=\left(l_{0} u\right)\left(x\right)+\int\limits _{t_{0} }^{t}\left(l_{1} u\right)\left(\tau \right)\, \beta _{1} \left(x\right)d\tau  +\int\limits _{t_{0} }^{t}\left(l_{2} u\right)\left(\tau \right)\, \beta _{2} \left(x\right)d\tau  +
$$
$$
+\int\limits _{t_{0} }^{t}\left(l_{3} u\right)\left(\tau \right)\, \beta _{3} \left(x\right)d\tau  +\iint \nolimits _{G}u_{txxx} \left(\tau ,\, \xi \right)R_{0} \left(\tau ,\, \xi ;\, t,x\right)d\tau\, d\xi  .
\end{equation}

Таким  образом, доказана  следующая.

\textbf{Теорема.} {\it Если $\det \gamma \ne 0$, то любая функция $u\in W_{p,n}^{\left(1,3\right)} \left(G\right)$ представима в виде (20).}

Очевидно, что правая часть формулы (20) определяется посредством значений $\left(l_{0} u\right)\left(x\right),\, \, \left(l_{1} u\right)\left(t\right),\, \, \left(l_{2} u\right)\left(t\right),$ $ \left(l_{3} u\right)\left(t\right)$ и  $u_{txxx} \left(t,x\right)$ операторов\linebreak $l_{0} ,\, l_{1} ,\, l_{2} ,\, l_{3} $ и $D_{t} D_{x}^{3} $ на рассматриваемой функции $u(t,x)$. Поэтому операторы $l_0,l_1,l_2,l_3$ и $D_tD_x^3$ будем называть  определяющими операторами для  представления (20).\\

{\bf 3. Задача Гурса классического вида как частный случай нелокальной комбинированной задачи (1), (2), (5).}

Теперь рассмотрим  следующий частный случай задачи (1), (2), (5):
\begin{equation} \label{GrindEQ__21_}
D_{t} D_{x}^{3} u\left(t,x\right)+\underset{i+j<4}{\sum\limits _{i=0}^{1}\sum\limits _{j=0}^{3}}
\left(D_{t}^{i} D_{x}^{j} u\left(t,x\right)\right)A_{i,j} \left(t,x\right)=  Z_{1,3} \left(t,x\right),\, \, \left(t,x\right)\in G
\end{equation}
\begin{equation} \label{GrindEQ__22_}
\left(l_{0} u\right)\left(x\right)\equiv u\left(t_{0} ,\, x\right)=Z_{0} \left(x\right),
\end{equation}
\begin{equation} \label{GrindEQ__23_}
\left\{\begin{array}{l} {\left(l_{1} u\right)\left(t\right)\equiv u_t\left(t,x_{0} \right)\alpha _{1,1} +u_{tx}\left(t, x_{0} \right)\alpha _{1,2} +u_{txx}\left(t, x_{0} \right)\alpha _{1,3} =Z_{1} \left(t\right)} \\
{\left(l_{2} u\right)\left(t\right)\equiv u_t\left(t,x_{0} \right)\alpha _{2,1} +u_{tx}\left(t, x_{0} \right)\alpha _{2,2} +u_{txx}\left(t, x_{0} \right)\alpha _{2,3} =Z_{2} \left(t\right)} \\
{\left(l_{3} u\right)\left(t\right)\equiv u_t\left(t,x_{0} \right)\alpha _{3,1} +u_{tx}\left(t, x_{0} \right)\alpha _{3,2} +u_{txx}\left(t, x_{0} \right)\alpha _{3,3} =Z_{3} \left(t\right)} \end{array}
\right.
\end{equation}

Пуст детерминант $3n$-мерной матрицы
\[
\gamma =\left(\begin{array}{l} {\alpha _{1,1} \, \, \, \, \alpha _{2,1} \, \, \, \alpha _{3,1} } \\ {\alpha _{1,2} \, \, \, \alpha _{2,2} \, \, \, \alpha _{3,2} } \\ {\alpha _{1,3} \, \, \, \, \alpha _{2,3} \, \, \, \, \alpha _{3,3} } \end{array}\right)
\]
отличен от нуля, т.е. $\det \gamma \ne 0$. Тогда  условия (23) можно привести  к виду
\begin{equation} \label{GrindEQ__24_}
\left\{\begin{array}{l} {u_t\left(t, x_{0} \right)=\bar{Z}_{1} \left(t\right)} \\ {u_{tx}\left(t,\, x_{0} \right)=\bar{Z}_{2} \left(t\right)} \\ {u_{txx}\left(t,\, x_{0} \right)=\bar{Z}_{3} \left(t\right)} \end{array}\right.
\end{equation}

Здесь:
\begin{equation} \label{GrindEQ__25_}
\left\{\begin{array}{l} {\bar{Z}_{1} \left(t\right)=Z_{1}(t) K_{1,1}^{0} +Z_{2} \left(t\right)K_{2,1}^{0} +Z_{3} \left(t\right)K_{3,1}^{0} } \\ {\bar{Z}_{2} \left(t\right)=Z_{1}(t) K_{1,2}^{0} +Z_{2} \left(t\right)K_{2,2}^{0} +Z_{3} \left(t\right)K_{3,2}^{0} } \\ {\bar{Z}_{3} \left(t\right)=Z_{1}(t) K_{1,3}^{0} +Z_{2} \left(t\right)K_{2,3}^{0} +Z_{3} \left(t\right)K_{3,3}^{0} } \end{array}\right.
\end{equation}
а
\[
K^{0} =\left(\begin{array}{l} {K_{1,1}^{0} \, \, \, \, K_{1,2}^{0} \, \, \, K_{1,3}^{0} } \\ {K_{2,1}^{0} \, \, \, \, K_{2,2}^{0} \, \, \, K_{2,3}^{0} } \\ {K_{3,1}^{0} \, \, \, \, K_{3,2}^{0} \, \, \, K_{3,3}^{0} } \end{array}\right)
\]
обратная матрица для матрицы $\gamma $, т.е. $K^{0} =\gamma ^{-1} $, причем  $K_{i,j}^{0} $- квадратные матрицы порядка $n$.

Равенства (25) показывают, что если $Z_{i} \in L_{p,n} \left(t_{0} ,\, t_{1} \right)$, $i=1,2,3$ тогда $\bar{Z}_{i} \in L_{p,n} \left(t_{0} ,\, t_{1} \right)$, $i=1,2,3$. Поэтому  условия (23) и (24) эквивалентны.

Кроме того, условия (22) и (24)  можно рассматривать как  условия  Гурса неклассического  вида. Очевидно также, что если $\det \gamma \ne 0$, то условия (22) и (23) эквиваленты условиям Гурса неклассического вида (22) и (24).

Если функция $u\in W_{p,n}^{\left(1,3\right)} \left(G\right)$ удовлетворяет условиям Гурса  неклассического вида (22) и (24),   то на удовлетворяет также  условиям Гурса  классического вида (22) и
\begin{equation} \label{GrindEQ__26_}
\left\{\begin{array}{l}
{u\left(t,\, x_{0} \right)=\varphi _{1} \left(t\right)} \\
{u_x\left(t,\, x_{0} \right)=\varphi _{2} \left(t\right)} \\
{u_{xx}\left(t,\, x_{0} \right)=\varphi _{3} \left(t\right)}
\end{array}\right.
\end{equation}
при
\begin{equation} \label{GrindEQ__27_}
\left\{\begin{array}{l} {\varphi _{1} \left(t\right)=Z_{0} \left(x_{0} \right)+\int\limits _{t_{0} }^{t}\bar{Z}_{1} \left(\tau \right)d\tau  } \\ {\varphi _{2} \left(t\right)=Z'_{0} \left(x_{0} \right)+\int\limits _{t_{0} }^{t}\bar{Z}_{2} \left(\tau \right)d\tau  } \\ {\varphi _{3} \left(t\right)=Z''_{0} \left(x_{0} \right)+\int\limits_{t_{0} }^{t}\bar{Z}_{3} \left(\tau \right)d\tau  } \end{array}\right.
\end{equation}

Верно и обратное, в том смысле,  что если функция $u\in W_{p,n}^{\left(1,3\right)} \left(G\right)$ удовлетворяет условиям (22) и (26) (где $Z_{0} \in W_{p,n}^{\left(3\right)} \left(x_{0} ,\, x_{1} \right)$ и $\varphi _{i} \in W_{p,n}^{\left(1\right)}
\left(t_{0},t_{1} \right),$ $i=1,2,3$) то  она удовлетворяет также  условиям (22), (24) при $\bar{Z}_{1} \left(t\right)=\varphi '_{1} \left(t\right)$, $\bar{Z}_{2} \left(t\right)=\varphi '_{2} \left(t\right),$ $\bar{Z}_{3} \left(t\right)=\varphi '_{3} \left(t\right)$. Таким образом, в пространстве $W_{p,n}^{\left(1,3\right)} \left(G\right)$ условия  Гурса неклассического вида (22), (24) эквиваленты условиям Гурса классического вида (22), (26). Однако, в случае условий Гурса (22), (26) правые части краевых условий, кроме условий

$Z_{0} \in W_{p,n}^{\left(3\right)} \left(x_{0} ,\, x_{1} \right)$ и $\varphi _{i} \in W_{p,n}^{\left(1\right)} \, \left(t_{0} ,\, t_{1} \right)$ должны подчиняться также  следующим условиям  согласования:\\
 $Z_{0} \left(x_{0} \right)=\varphi _{1} \left(t_{0} \right)$, $Z'_{0} \left(x_{0} \right)=\varphi _{2} \left(t_{0} \right)$, $Z''_{0} \left(x_{0} \right)=\varphi _{3} \left(t_{0} \right)$.

В случае же условий Гурса неклассического вида (22), (24)  на правых частей краевых условий, кроме $Z_{0} \in W_{p,n}^{\left(3\right)} \left(x_{0} ,\, x_{1} \right)$ и  $\bar{Z}_{i} \in L_{p,n} \left(t_{0} ,\, t_{1} \right),\, i=1,2,3$,  никаких дополнительных условий  типа согласования не требуются. Поэтому  задача Гурса неклассического вида (21), (22), (24)  по постановке является более естественной, чем задача Гурса  классического вида (21), (22), (26).

\newpage

\begin{titlepage}

\centerline{\bf Сведения об авторе}

\

{\bf 1. Фамилия, имя, отчество:} Мамедов Ильгар Гурбат оглы.

{\bf 2. Называние организации:} Институт Кибернетики им. А.И.Гусейнова НАН Азербайджана.

{\bf 3. Дольжность, уч. степень, уч. звание:} ведущий научный сотрудник, кандидат физ.-мат. наук, доцент.

{\bf 4. Почтовый адрес:} Азербайджан, AZ1141, г.Баку, ул. Б.Вагабзаде, д.9,
Институт Кибернетики им. А.И.Гусейнова НАН Азербайджана.

{\bf 5.	Телефон, E-mail, факс:}  (994 12) 539 28 26, \\
ilgar-mammadov@rambler.ru, (994 12) 539 28 26,

\end{titlepage}

\begin{thebibliography}{99}

\bibitem{1}
А.В.Бицадзе, А.А.Самарский,  О некоторых простейших обобщениях  линейных эллиптических краевых задач // ДАН СССР, т. 185, \No 4,  с. 739-740, 1969.

\bibitem{2}
Н.И.Ионкин,  Решение одной краевой задачи теории теплопроводности с неклассическим краевым условиям // ДУ, т.13, \No 2, с. 294-304, 1977.

\bibitem{3}
 Н.И.Ионкин, Е.И.Моисеев,  О задаче для уравнения  теплопроводности  с двуточечными  краевыми условиями // ДУ, т.15, \No 7, с. 1284-1296, 1979.

\bibitem{4}
А.А.Самарский,  О некоторых проблемах теории современной дифференциальных уравнений // ДУ, т. 16, \No 11, с. 1221-1228, 1980.

\bibitem{5}
 М.Х.Шхануков, А.П.Солдатов,   Краевые задачи с общим  нелокальным условием А.А. Самарского  для псевдопараболических уравнений  высокого порядка // ДАН СССР, т.297, \No 3, с.547-552, 1987.

\bibitem{6}
R.E.Showalter, T.W.Ting Pseudo-parabolic partial differential equations // Math. Analys. v.1, pp.1-26, 1970.


\bibitem{7}
D.Colton,  Pseudoparabolic equations in one space variable // J. Different. Equat.  v.12, No 3, pp. 559-565, 1972.

\bibitem{8}
W.Rundell, M.Stecher,  The uniqueness class for the Cauchy problem for pseudoparabolic equations // Proc. Amer. Math. Soc.  v.76, \No 2, pp. 253-257, 1979.


\bibitem{9}
В.Д.Гилев, Г.А.Шадрин,  Построение фундаментального решения для уравнения, описывающего движение жидкости в трещиноватых средах // Вычисл. матем. и программирование. М.: Изд-во МГПИ им. В.И.Ленина, Вып. 4, с.102-111, 1976.

\bibitem{10}
М.Х.Шхануков,  О некоторых  краевых задачах для  уравнения  третьего порядка, возникающих при моделировании фильтрации жидкости в  пористых средах // ДУ, т.18, \No 4, с. 689-699, 1982.

\bibitem{11}
М.Х.Шхануков,  О некоторых  краевых задачах для  уравнения  третьего порядка и экстремальных  свойствах его решений //  ДУ, т.19, \No 1, с. 145-152, 1983.

\bibitem{12}
Н.М.Салтыкова,  Обобщенная задача Бицадзе--Самарского для уравнения смешанного типа второго рода // Известия высших учебных заведений. Математика, No 11 (234), с. 43-48, 1981.

\bibitem{13}
Ф.Ш.Ахмедов,   Оптимизация гиперболических систем при  неклокальных краевых условиях типа  Бицадзе--Самарского // ДАН СССР, т.283, \No 4, с. 787-791, 1985.

\bibitem{14}
 О.А.Репин, А.А.Килбас,  Аналог задачи Бицадзе--Самарского для уравнения смешанного типа с дробной производной // Дифференциальные уравнения, т.39, \No 5, с. 638-644, 2003.

\bibitem{15}
Л.А.Ковалева,  О модифицированной задаче Бицадзе-Самарского // Вестник Самарского Государственного Технического Университета. Сер.: физико-математические науки, \No 1, с. 10-15, 2007.

\bibitem{16}
К.Б.Сабитов,  Краевая задача для уравнения параболо-гиперболического типа с нелокальным интегральным условием // Дифференциальные уравнения, т. 46, \No 10, с.1468-1478, 2010.


\bibitem{17}
А.В.Гулин, Н.С.Удовиченко,  Разностная схема для задачи Самарского-Ионкина с параметром // Дифференциальные уравнения, т. 44, \No 7, с.963-969, 2008.


\bibitem{18}
I.G.Mamedov,  Generalization of  multipoint boundary-value problems  of Bitsadze--Samarski and Samarski--Ionkin type for fourth order  loaded hyperbolic integro-differential equations and their  operator generalization // Proc. of IMM of NAS of Azerbaijan, v. XXIII, pp. 77-84, 2005.

\bibitem{19}
И.Г.Мамедов,  Смешанная задача с нелокальными краевыми  условиями типа Бицадзе--Самарского и Самарского--Ионкина,  возникающая  при моделировании фильтрации жидкости в трещиноватых средах // Известия НАН Азерб., сер. физ.--техн. и матем. наук., т. XXVI, \No 3, С.32-37, 2006.


\bibitem{20}
Ю.М.Березанский,  Я.А.Ройтберг,  {  Теорема о гомеоморфизмах и функция Грина для общих эллиптических граничных задач} // {Укр. мат. жур.}, т. 19, \No 5, с.3-32, 1967.

\bibitem{21}
Н.В.Житарашу,  {  Теорема о полном наборе изоморфизмов в $L_2$-теории модельных начальных параболических краевых задач} // Математические исследования, {Кишинев}.   \No 88, с. 40-59, 1986.


\bibitem{22}
С.С.Ахиев,  Фундаментальные решения некоторых локальных  и нелокальных краевых задач  и их представления // ДАН СССР, т.271, \No 2, с. 265-269, 1983.

\bibitem{23}
{И.Г.Мамедов }
Фундаментальное решение задачи Коши, связанной с псевдопараболическим уравнением четвертого порядка // ЖВМ и МФ, т. 49, \No 1, с. 99-110, 2009.

\bibitem{24}
И.Г.Мамедов,  Об одной задаче Гурса в пространстве Соболева // Известия вузов. Математика.  \No 2, с.54-64, 2011.
\end{thebibliography}
\end{document}